\documentclass[11pt]{article}
\usepackage{amssymb}

\usepackage{amsmath}
\usepackage{theorem}
\usepackage{here}
\usepackage[dvipdfmx]{color}

\newtheorem{thm}{Theorem}[section]
\newtheorem{lem}{Lemma}[section]

\theorembodyfont{\rmfamily}

\makeatletter
\renewcommand{\theequation}{%
   \thesection.\arabic{equation}}
\@addtoreset{equation}{section}
\makeatother

\usepackage{comment} %
\usepackage{bm}
\usepackage[dvipdfmx]{graphicx} %

\def\be{{\beta}}
\def\ga{{\gamma}}
\def\de{{\delta}}
\def\ep{{\varepsilon}}

\def\si{{\sigma}}
\def\om{{\omega}}

\def\bbe{{\text{\boldmath $\beta$}}}

\def\Ga{{\Gamma}}

\def\a{{\text{\boldmath $a$}}}
\def\b{{\text{\boldmath $b$}}}

\def\e{{\text{\boldmath $e$}}}

\def\i{{\text{\boldmath $i$}}}

\def\v{{\text{\boldmath $v$}}}

\def\x{{\text{\boldmath $x$}}}
\def\y{{\text{\boldmath $y$}}}
\def\z{{\text{\boldmath $z$}}}

\def\X{{\text{\boldmath $X$}}}

\def\Z{{\text{\boldmath $Z$}}}

\def\vbt{{\tilde \v}}

\def\abt{{\tilde \a}}
\def\bbt{{\tilde \b}}

\def\Ih{{\widehat I}}

\def\Xbt{{\widetilde \X}}

\def\Zbt{{\widetilde \Z}}

\def\Ic{{\cal I}}

\def\Lc{{\cal L}}
\def\Kc{{\cal K}}
\def\[{{\text{\boldmath $[$}}}
\def\]{{\text{\boldmath $]$}}}

\def\/{{\Bigr/\!\!}}

\def\1r{{\rm (1)}}
\def\2r{{\rm (2)}}
\def\3r{{\rm (3)}}
\def\4r{{\rm (4)}}
\def\5r{{\rm (5)}}

\def\non{{\nonumber}}

\begin{document}
\title{A Simple Proof of Posterior Robustness}
\author{
Yasuyuki Hamura\footnote{Graduate School of Economics, Kyoto University, 
Yoshida-Honmachi, Sakyo-ku, Kyoto, 606-8501, JAPAN. 
\newline{
E-Mail: yasu.stat@gmail.com}} \
}
\maketitle
\begin{abstract}
Conditions for Bayesian posterior robustness have been examined in recent literature. 
However, many of the proofs seem to be long and complicated. 
In this paper, we first summarize some basic lemmas that have been applied implicitly or explicitly. 
Then, using them, we give a simple proof of posterior robustness. 
Our sufficient condition is new and practically relevant.

\par\vspace{4mm}
{\it Key words and phrases:\ Bayes, linear regression, posterior robustness, super heavy-tailed distribution. } 
\end{abstract}

\section{Introduction}
\label{sec:introduction}
In applications of Bayesian linear regression models, we often want to base our posterior inference on moderate observations after discarding outlying part of the data, which is regarded as not directly related to parameters of interest. 
Since it is not necessarily clear to us which observations are extreme, it is desirable that the effects of outliers be automatically removed in a Bayesian way. 
Research on such Bayesian posterior robustness (O'Hagan, 1979) and related topics has a long history; see, for example, West (1984), Andrade and O'Hagan (2006, 2011), and O'Hagan and Pericchi (2012). 
However, there have been relatively few theoretical studies until the fundamental and pathbreaking papers by Desgagn\'{e} (2013, 2015) and Gagnon et al. (2019). 

Desgagn\'{e} (2013, 2015) first introduced log-regularly varying distributions, which are crucial for robust Bayesian analysis, and derived conditions for (full) posterior robustness for a simple (location-)scale model. 
Then Gagnon et al. (2019) obtained results for a general regression model. 
Since then, much work has been done to establish posterior robustness in a variety of important settings. 
For example, Hamura et al. (2022) considered using novel error distributions and shrinkage priors on regression coefficients, while Hamura et al. (2021) examined the discrete case. 
Andrade (2022) and Gagnon and Hayashi (2023) considered partially robust Student's $t$-models. 
However, it seems that many of the proofs %
are long and complicated and not easy to follow for some researchers. 
Also, available conditions for posterior robustness can be restrictive in some cases. 

The main purpose of this paper is to summarize some basic results that have been applied implicitly or explicitly in the literature and give a simple proof of posterior robustness to illustrate their use. 
This could make the topic more accessible. 
Additionally, the sufficient condition we derive is a new one and can be useful when the proportion of outliers may be large if a heavy-tailed prior is assumed for a scale parameter. 
This is another contribution of the paper. 

The remainder of the paper is organized as follows. 
In Section \ref{sec:model}, the Bayesian regression model we consider is described and a definition of posterior robustness is given. 
In Section \ref{sec:conditions}, our sufficient condition for posterior robustness is presented and compared with conditions avaliable in the literature. 
In Section \ref{sec:proof}, two key results are explicitly stated and the main theorem of Section \ref{sec:conditions} is proved. 
In the Appendix, a detailed proof of a lemma is given. 
Additional results used in Section \ref{sec:proof} are given in the Supplementary Material.

\section{The Model}
\label{sec:model} 
Suppose that for $i = 1, \dots , n$, we observe 
\begin{align}
&y_i \sim f(( y_i - {\x _i}^{\top } \bbe ) / \si ) / \si \text{,} \non 
\end{align}
where $\x _i$ are continuous explanatory variables while $\bbe \in \mathbb{R} ^p$ and $\si \in (0, \infty )$ are regression coefficients and a scale parameter, respectively, and where $f$ is a (proper) error density. 
Let our prior distribution for $( \bbe , \si )$ be denoted by 
\begin{align}
( \bbe , \si ) \sim \pi ( \bbe , \si ) \text{.} \non 
\end{align}
Then, for any sets of observations $\y _{\Ic } = \{ y_i | i \in \Ic \} \subset \y = \{ y_1 , \dots , y_n \} $, where $\Ic \subset I = \{ 1, \dots , n \} $, the posterior density of $( \bbe , \si )$ given $\y _{\Ic }$ is 
\begin{align}
p( \bbe , \si | \y _{\Ic } ) &= \frac{ \pi ( \bbe , \si ) \prod_{i \in \Ic } \{ f(( y_i - {\x _i}^{\top } \bbe ) / \si ) / \si \} }{ \int_{\mathbb{R} ^p \times (0, \infty )} \pi ( \bbe , \si ) \big[ \prod_{i \in \Ic } \{ f(( y_i - {\x _i}^{\top } \bbe ) / \si ) / \si \} \big] d( \bbe , \si ) } \text{.} \non 
\end{align}

Following Desgagn\'{e} (2015) and other papers, let $a_1 , \dots , a_n \in \mathbb{R}$ and $b_1 , \dots , b_n \in \mathbb{R}$ and suppose that for all $i = 1, \dots , n$, we have $y_i = a_i + b_i \om $ and $\om \to \infty $. 
Let $\Kc , \Lc \subset I$ be such that $\Kc \neq \emptyset $, $\Lc \neq \emptyset $, $\Kc \cap \Lc = \emptyset $, and $\Kc \cup \Lc = I$ and suppose that $\Lc = \{ i \in I | b_i \neq 0 \} $. 
Thus, $\Kc $ and $\Lc $ are interpreted as the sets of indices for nonoutlying and outlying observations, respectively. 

Note that $\y $ depends on $\om \to \infty $ while $\y _{\Kc }$ does not. 
We say that the posterior is robust to the outliers if 
\begin{align}
p( \bbe , \si | \y ) \to p( \bbe , \si | \y _{\Kc } ) \non 
\end{align}
as $\om \to \infty $ at each $( \bbe , \si ) \in \mathbb{R} ^p \times (0, \infty )$. 
Here, $\Kc $ and $\Lc $ %
are assumed to be unknown. 
Thus, the posterior robustness means that the effects of the outliers are automatically removed from our posterior inference based on $p( \bbe , \si | \y )$. 

The question is under what conditions on $f \colon \mathbb{R} \to (0, \infty )$, $\pi \colon \mathbb{R} ^p \times (0, \infty ) \to (0, \infty )$, and $| \Lc | = n - | \Kc |$ the posterior robustness holds so that we have 
\begin{align}
p( \bbe , \si | \y ) &= \frac{ \displaystyle \pi ( \bbe , \si ) \Big[ \prod_{i \in \Kc } \{ f(( y_i - {\x _i}^{\top } \bbe ) / \si ) / \si \} \Big] \prod_{i \in \Lc } {f(( y_i - {\x _i}^{\top } \bbe ) / \si ) / \si \over f( y_i )} }{ \displaystyle \int_{\mathbb{R} ^p \times (0, \infty )} \pi ( \bbe , \si ) \Big[ \prod_{i \in \Kc } \{ f(( y_i - {\x _i}^{\top } \bbe ) / \si ) / \si \} \Big] \Big\{ \prod_{i \in \Lc } {f(( y_i - {\x _i}^{\top } \bbe ) / \si ) / \si \over f( y_i )} \Big\} d( \bbe , \si ) } \non \\
&\to \frac{ \displaystyle \pi ( \bbe , \si ) \prod_{i \in \Kc } \{ f(( y_i - {\x _i}^{\top } \bbe ) / \si ) / \si \} }{ \displaystyle \int_{\mathbb{R} ^p \times (0, \infty )} \pi ( \bbe , \si ) \Big[ \prod_{i \in \Kc } \{ f(( y_i - {\x _i}^{\top } \bbe ) / \si ) / \si \} \Big] d( \bbe , \si ) } = p( \bbe , \si | \y _{\Kc } ) \non 
\end{align}
as $\om \to \infty $. 
Informally, it has turned out in the literature that with regard to $f$, it is necessary and sufficient that it is log-regularly varying / super heavy-tailed or that we have $f(|z|) \approx |z|^{- 1}$ as $|z| \to \infty $ ignoring log factors. 
This restriction is to ensure that $\{ f(( y_i - {\x _i}^{\top } \bbe ) / \si ) / \si \} f( y_i ) \to 1$ for all $i \in \Lc $ so that the numerator of the above expression converges to the correct (unnormalized) density. 
In contrast, available conditions on $\pi $ and $| \Lc |$ tend to be dependent. 
They are imposed to justify the interchange of limit and integral in the denominator. %

\section{Conditions for Posterior Robustness}
\label{sec:conditions} 
In what follows, we assume for simplicity that 
\begin{align}
f(z) = {\ga / 2 \over 1 + |z|} {1 \over \{ 1 + \log (1 + |z|) \} ^{1 + \ga }} \non 
\end{align}
for all $z \in \mathbb{R}$ for some $\ga > 0$.

\begin{thm}
\label{thm:proper_bounded} 
Suppose that %
$\int_{\mathbb{R} ^p \times (0, \infty )} (1 + \si ^{\rho } ) \pi ( \bbe , \si ) d( \bbe , \si ) < \infty $ for some $\rho > 0$. 
Suppose that there exist $\nu > 0$ and $M > 0$ such that 
\begin{align}
p( \bbe | \si ) &\le M \prod_{k = 1}^{p} \Big\{ {1 \over \si } {1 \over (1 + | \be _k | / \si )^{1 + \nu }} \Big\} \label{eq:assumption_conditional_prior} %
\end{align}
for all $\bbe = ( \be _1 , \dots , \be _p ) \in \mathbb{R} ^p$ and all $\si \in (0, \infty )$. 
Suppose that $| \Kc | \ge | \Lc | + p$. 
Then we have 
\begin{align}
\lim_{\om \to \infty } p( \bbe , \si | \y ) &= p( \bbe , \si | \y _{\Kc } ) \non 
\end{align}
at each $( \bbe , \si ) \in \mathbb{R} ^p \times (0, \infty )$. 
\end{thm}

The first condition is satisfied for most proper priors. 
Condition (\ref{eq:assumption_conditional_prior}) is satisfied, for example, if we use a conditionally independent prior such that 
\begin{align}
p( \be _k | \si ) &\propto {1 \over \si } {1 \over (1 + | \be _k | / \si )^{1 + \nu _k}} \label{eq:prior_independent_t} 
\end{align}
for some $\nu _k > 0$ for all $k = 1, \dots , p$. 
The condition is satisfied also when $\bbe | \si $ is multivariate $t$ with 
\begin{align}
p( \bbe | \si ) &= {1 \over \si ^p} {\Ga (( \nu + p) / 2) \over \Ga ( \nu / 2) \nu ^{p / 2} \pi ^{p / 2}} {1 \over (1 + \| \bbe / \si \| ^2 / \nu )^{( \nu + p) / 2}} \non \\
&\le M \prod_{k = 1}^{p} \Big\{ {1 \over \si } {1 \over (1 + | \be _k | / \si )^{1 + \nu / p}} \Big\} \non 
\end{align}
for all $\bbe = ( \be _1 , \dots , \be _p ) \in \mathbb{R} ^p$ and all $\si \in (0, \infty )$ for some $\nu , M > 0$. 

Although Theorem 2.1 of Gagnon et al. (2019) is applicable also to improper priors, they assume that $| \Kc | \ge | \Lc | + 2 p - 1$, which is more restrictive than the condition on $| \Lc |$ given above. 
Additionally, Gagnon et al. (2019) assume that $\pi ( \bbe , \si ) \le M' \max \{ 1, 1 / \si \} $ for all $( \bbe , \si ) \in \mathbb{R} ^p \times (0, \infty )$ for some $M' > 0$, which is not satisfied for some choice of the prior for $\si $ when we use (\ref{eq:prior_independent_t}), for example. 
On the other hand, Hamura et al. (2022) assume the same condition on $| \Lc |$ as that of Theorem \ref{thm:proper_bounded} and a different condition on the prior distribution of $\bbe | \si $ but their main result is applicable only if $\si $ has sufficiently many negative moments. 
Roughly speaking, Hamura et al. (2022) consider unbounded prior densites of $\bbe | \si $, whereas we consider unbounded prior densites of $\si $. 
Thus, Theorem \ref{thm:proper_bounded} is a new result that is relevant to practical situations.

\section{Proof of Theorem \ref{thm:proper_bounded}}
\label{sec:proof} 
The following result is due to Gagnon et al. (2019) but proved in the Appendix for clarity. 

\begin{lem}
\label{lem:gagnon} 
Let $m, p \in \mathbb{N}$ satisfy $m \ge p + 1$. 
Let $\z _1 , \dots , \z _m \in \mathbb{R} ^p$ and $( a_1 , b_1 ), \dots , ( a_m , b_m ) \in \mathbb{R} ^2$ satisfy the following three conditions: 
\begin{itemize}
\item[{\rm{(i)}}]
For any set of $p$ distinct indices $i_1 , \dots , i_p \in \{ 1, \dots , m \} $, $( \z _{i_1} , \dots , \z _{i_p} )^{\top }$ is invertible. 
\item[{\rm{(ii)}}]
For any set of $p + 1$ distinct indices $i_1 , \dots , i_{p + 1} \in \{ 1, \dots , m \} $, $(( \z _{i_1} , \dots , \z _{i_{p + 1}} )^{\top } , ( a_{i_1} , \dots , a_{i_{p + 1}} )^{\top } )$ is invertible. 
\item[{\rm{(iii)}}]
For any set of $p + 1$ distinct indices $i_1 , \dots , i_{p + 1} \in \{ 1, \dots , m \} $, $(( \z _{i_1} , \dots , \z _{i_{p + 1}} )^{\top } , ( b_{i_1} , \dots , b_{i_{p + 1}} )^{\top } )$ is invertible unless $b_{i_1} = \dots = b_{i_{p + 1}} = 0$. 
\end{itemize}
Then there exist $\ep > 0$ and $M > 0$ such that for all $\om \ge M$, we have 
\begin{align}
\mathbb{R} ^p \subset \bigcup_{i_1 = 1}^{m} \dots \bigcup_{i_p = 1}^{m} \bigcap_{\{ 1, \dots , m \} \setminus \{ i_1 , \dots , i_p \} } \{ \bbe \in \mathbb{R} ^p | | a_i + b_i \om - {\z _i}^{\top } \bbe | > \ep \} \text{.} \non 
\end{align}
\end{lem}

The following lemma is exactly as given in the Supplementary material of Hamura et al. (2022). 

\begin{lem}
\label{lem:linear} 
Let $m, p \in \mathbb{N}$. 
Let %
$( w_1 , \dots , w_m )^{\top } \in \mathbb{R} ^m$. 
Let %
$( \z _1 , \dots , \z _m )^{\top } \in \mathbb{R} ^{m \times p}$ be such that any set of its $p$ distinct row vectors is linearly independent. 
Suppose that $m \ge p$. 
Then there exist $R > 0$ and $\de > 0$ %
such that 
\begin{align}
\prod_{i = 1}^{m} {1 \over 1 + | w_i - {\z _i}^{\top } \bbe |} \le {1 \over (1 + \de \| \bbe \| )^{m - p + 1}} \non 
\end{align}
for all $\bbe \in \mathbb{R} ^p$ satisfying $\| \bbe \| \ge R$. 
\end{lem}

We now prove Theorem \ref{thm:proper_bounded}. 

\bigskip

\noindent
{\bf Proof of Theorem \ref{thm:proper_bounded}.} \ \ By part (i) of Lemma \ref{lem:f}, 
\begin{align}
\lim_{\om \to \infty } {f(( y_i - {\x _i}^{\top } \bbe ) / \si ) / \si \over f( y_i )} = 1 \non 
\end{align}
for all $i \in \Lc $. 
Therefore, it is sufficient to show that 
\begin{align}
\lim_{\om \to \infty } \int_{\mathbb{R} ^p \times (0, \infty )} h( \bbe , \si ; \om ) d( \bbe , \si ) = \int_{\mathbb{R} ^p \times (0, \infty )} \pi ( \bbe , \si ) \Big[ \prod_{i \in \Kc } \{ f(( y_i - {\x _i}^{\top } \bbe ) / \si ) / \si \} \Big] d( \bbe , \si ) \text{,} \non 
\end{align}
where 
\begin{align}
h( \bbe , \si ; \om ) &= \pi ( \bbe , \si ) \Big[ \prod_{i \in \Kc } \{ f(( y_i - {\x _i}^{\top } \bbe ) / \si ) / \si \} \Big] \prod_{i \in \Lc } {f(( y_i - {\x _i}^{\top } \bbe ) / \si ) / \si \over f( y_i )} \non \\
&= \Big\{ \prod_{i \in \Lc } 1(| y_i - {\x _i}^{\top } \bbe | \ge | y_i | / 2) \Big\} h( \bbe , \si ; \om ) + \Big\{ 1 - \prod_{i \in \Lc } 1(| y_i - {\x _i}^{\top } \bbe | \ge | y_i | / 2) \Big\} h( \bbe , \si ; \om ) \non 
\end{align}
for $( \bbe , \si ) \in \mathbb{R} ^p \times (0, \infty )$ for $\om > 0$. 

By part (ii) of Lemma \ref{lem:f}, 
\begin{align}
&\Big\{ \prod_{i \in \Lc } 1(| y_i - {\x _i}^{\top } \bbe | \ge | y_i | / 2) \Big\} h( \bbe , \si ; \om ) \non \\
&\le \Big\{ \prod_{i \in \Lc } 1(| y_i - {\x _i}^{\top } \bbe | \ge | y_i | / 2) \Big\} M_1 \pi ( \bbe , \si ) \{ 1 + \log (1 + \si ) \} ^{| \Lc | (1 + \ga )} \prod_{i \in \Kc } \{ f(( y_i - {\x _i}^{\top } \bbe ) / \si ) / \si \} \non \\
&\le M_2 \pi ( \bbe , \si ) (1 + \si ^{\rho } ) \prod_{i \in \Kc } \{ f(( y_i - {\x _i}^{\top } \bbe ) / \si ) / \si \} \non 
\end{align}
for sufficiently large $\om > 0$ and all $( \bbe , \si ) \in \mathbb{R} ^p \times (0, \infty )$ for some $M_1 , M_2 > 0$. 
Therefore, by the dominated convergence theorem, 
\begin{align}
&\lim_{\om \to \infty } \int_{\mathbb{R} ^p \times (0, \infty )} \Big\{ \prod_{i \in \Lc } 1(| y_i - {\x _i}^{\top } \bbe | \ge | y_i | / 2) \Big\} h( \bbe , \si ; \om ) d( \bbe , \si ) \non \\
&= \int_{\mathbb{R} ^p \times (0, \infty )} \Big[ \lim_{\om \to \infty } \Big\{ \prod_{i \in \Lc } 1(| y_i - {\x _i}^{\top } \bbe | \ge | y_i | / 2) \Big\} h( \bbe , \si ; \om ) \Big] d( \bbe , \si ) \non \\
&= \int_{\mathbb{R} ^p \times (0, \infty )} \pi ( \bbe , \si ) \Big[ \prod_{i \in \Kc } \{ f(( y_i - {\x _i}^{\top } \bbe ) / \si ) / \si \} \Big] d( \bbe , \si ) \text{.} \non 
\end{align}
Thus, it suffices to prove that 
\begin{align}
&\lim_{\om \to \infty } \int_{\mathbb{R} ^p \times (0, \infty )} \Big\{ 1 - \prod_{i \in \Lc } 1(| y_i - {\x _i}^{\top } \bbe | \ge | y_i | / 2) \Big\} h( \bbe , \si ; \om ) d( \bbe , \si ) = 0 \text{.} \non 
\end{align}

Let $\Ih = \{ - 1, \dots , - p \} \cup \{ 1, \dots , n \} $. 
Let $( a_i , b_i ) = (0, 0)$ and $\x _i = \e _{- i}^{(p)}$ for $i \in \{ - 1, \dots , - p \} $. 
Then, by Lemma \ref{lem:gagnon}, there exist $\ep , M > 0$ such that 
\begin{align}
\mathbb{R} ^p \subset %
\bigcup_{\i %
\in {\Ih {}}^p} A_{\i } ( \om ) \non 
\end{align}
for all $\om \ge M$, where 
\begin{align}
A_{( i_1 , \dots , i_p )} ( \om ) &= \bigcap_{i \in \Ih \setminus \{ i_1 , \dots , i_p \} } \{ \bbe \in \mathbb{R} ^p \mid | a_i + b_i \om - {\x _i}^{\top } \bbe | > \ep \} \non 
\end{align}
for $\om > 0$ for $( i_1 , \dots , i_p ) \in {\Ih {}}^p$. 
Therefore, it is sufficient to show that 
\begin{align}
&\lim_{\om \to \infty } \int_{\mathbb{R} ^p \times (0, \infty )} h_{\i } ( \bbe , \si ; \om ) d( \bbe , \si ) = 0 \non 
\end{align}
for all $\i \in {\Ih {}}^p$, 
where 
\begin{align}
h_{\i } ( \bbe , \si ; \om ) &= 1( \bbe \in A_{\i } ( \om )) \Big\{ 1 - \prod_{i \in \Lc } 1(| y_i - {\x _i}^{\top } \bbe | \ge | y_i | / 2) \Big\} h( \bbe , \si ; \om ) \non 
\end{align}
for $( \bbe , \si ) \in \mathbb{R} ^p \times (0, \infty )$ for $\om > 0$ for $\i \in {\Ih {}}^p$. 

Fix $\i = ( i_1 , \dots , i_p ) \in {\Ih {}}^p$. 
Let $- p \le i(1) < \dots < i(p) \le n$ be such that $\{ i_1 , \dots , i_p \} \subset \{ i(1), \dots , i(p) \} \subset \Ih $. 
Then 
\begin{align}
h_{\i } ( \bbe , \si ; \om ) &\le 1( \bbe \in A_{\i } ( \om )) \Big\{ 1 - \prod_{i \in \Lc } 1(| y_i - {\x _i}^{\top } \bbe | \ge | y_i | / 2) \Big\} \non \\
&\quad \times M \pi ( \si ) \prod_{k = 1}^{p} \Big\{ {1 \over \si } {1 \over (1 + | \be _k | / \si )^{1 + \nu }} \Big\} \Big[ \prod_{i \in \Kc } \{ f(( y_i - {\x _i}^{\top } \bbe ) / \si ) / \si \} \Big] \prod_{i \in \Lc } {f(( y_i - {\x _i}^{\top } \bbe ) / \si ) / \si \over f( y_i )} \non \\
&= 1( \bbe \in A_{\i } ( \om )) \Big\{ 1 - \prod_{i \in \Lc } 1(| y_i - {\x _i}^{\top } \bbe | \ge | y_i | / 2) \Big\} M \pi ( \si ) \prod_{j = 1}^{p} g_j ( \bbe , \si ; \om ) \non \\
&\quad \times \Big[ \prod_{i \in \{ - 1, \dots , - p \} \setminus \{ i(1), \dots , i(p) \} } \Big\{ {1 \over \si } {1 \over (1 + | \be _{- i} | / \si )^{1 + \nu }} \Big\} \Big] \non \\
&\quad \times \Big[ \prod_{i \in \Kc \setminus \{ i(1), \dots , i(p) \} } \{ f(( y_i - {\x _i}^{\top } \bbe ) / \si ) / \si \} \Big] {\prod_{i \in \Lc \setminus \{ i(1), \dots , i(p) \} } f(( y_i - {\x _i}^{\top } \bbe ) / \si ) / \si \over \prod_{i \in \Lc } f( y_i )} \text{,} \non 
\end{align}
where 
\begin{align}
g_j ( \bbe , \si ; \om ) &= \begin{cases} \displaystyle f( \{ y_{i(j)} - {\x _{i(j)}}^{\top } \bbe \} / \si ) / \si \text{,} & \text{if $i(j) \ge 1$} \text{,} \\ \displaystyle {1 \over \si } {1 \over \{ 1 + | \be _{- i(j)} | / \si \} ^{1 + \nu }} \text{,} & \text{if $i(j) \le - 1$} \text{,} \end{cases} \non 
\end{align}
for $( \bbe , \si ) \in \mathbb{R} ^p \times (0, \infty )$ for $\om > 0$ for $j = 1, \dots , p$. 
For any $i \in \{ - 1, \dots , - p \} \setminus \{ i(1), \dots , i(p) \} $, 
\begin{align}
1( \bbe \in A_{\i } ( \om )) {1 \over \si } {1 \over (1 + | \be _{- i} | / \si )^{1 + \nu }} &\le 1( \bbe \in A_{\i } ( \om )) {1 \over \si } {1 \over 1 + | a_i + b_i \om - {\x _i}^{\top } \bbe | / \si } \non \\
&\le 1( \bbe \in A_{\i } ( \om )) {1 \over \si } {1 \over ( \ep + | a_i + b_i \om - {\x _i}^{\top } \bbe |) / (2 \si )} \non \\
&= 1( \bbe \in A_{\i } ( \om )) {2 \over \ep + |0 - ( \e _{- i}^{(p)} )^{\top } \bbe |} \non 
\end{align}
for all $\om > 0$ and all $( \bbe , \si ) \in \mathbb{R} ^p \times (0, \infty )$. 
For any $i \in \Kc \setminus \{ i(1), \dots , i(p) \} $, by part (iii) of Lemma \ref{lem:f}, 
\begin{align}
1( \bbe \in A_{\i } ( \om )) f(( y_i - {\x _i}^{\top } \bbe ) / \si ) / \si &\le 1( \bbe \in A_{\i } ( \om )) M_3 {1 \over \si } {1 \over 1 + | y_i - {\x _i}^{\top } \bbe | / \si } \non \\
&\le 1( \bbe \in A_{\i } ( \om )) M_3 {1 \over \si } {1 \over ( \ep + | y_i - {\x _i}^{\top } \bbe |) / (2 \si )} \non \\
&= 1( \bbe \in A_{\i } ( \om )) M_3 {2 \over \ep + | a_i - {\x _i}^{\top } \bbe |} \non 
\end{align}
for all $\om > 0$ and all $( \bbe , \si ) \in \mathbb{R} ^p \times (0, \infty )$, where $M_3 > \ga / 2$. 
For any $i \in \Lc \setminus \{ i(1), \dots , i(p) \} $, by part (iii) of Lemma \ref{lem:f}, 
\begin{align}
1( \bbe \in A_{\i } ( \om )) f(( y_i - {\x _i}^{\top } \bbe ) / \si ) / \si &\le 1( \bbe \in A_{\i } ( \om )) M_3 {1 \over \si } {1 \over 1 + | y_i - {\x _i}^{\top } \bbe | / \si } \non \\
&\le 1( \bbe \in A_{\i } ( \om )) M_3 {1 \over \si } {1 \over ( \ep + | y_i - {\x _i}^{\top } \bbe |) / (2 \si )} \non \\
&\le 1( \bbe \in A_{\i } ( \om )) M_3 {2 \over \ep } \non 
\end{align}
for all $\om > 0$ and all $( \bbe , \si ) \in \mathbb{R} ^p \times (0, \infty )$. 
Therefore, 
\begin{align}
&h_{\i } ( \bbe , \si ; \om ) \non \\
&\le 1( \bbe \in A_{\i } ( \om )) \Big\{ 1 - \prod_{i \in \Lc } 1(| y_i - {\x _i}^{\top } \bbe | \ge | y_i | / 2) \Big\} \pi ( \si ) \prod_{j = 1}^{p} g_j ( \bbe , \si ; \om ) \non \\
&\quad \times M_4 \Big\{ \prod_{i \in \{ - 1, \dots , - p \} \setminus \{ i(1), \dots , i(p) \} } {1 \over \ep + |0 - ( \e _{- i}^{(p)} )^{\top } \bbe |} \Big\} \Big\{ \prod_{i \in \Kc \setminus \{ i(1), \dots , i(p) \} } {1 \over \ep + | a_i - {\x _i}^{\top } \bbe |} \Big\} {1 \over \prod_{i \in \Lc } f( y_i )} \non 
\end{align}
for all $\om > 0$ and all $( \bbe , \si ) \in \mathbb{R} ^p \times (0, \infty )$ for some $M_4 > 0$. 
Note that by Lemma \ref{lem:linear}, 
\begin{align}
&\Big\{ \prod_{i \in \{ - 1, \dots , - p \} \setminus \{ i(1), \dots , i(p) \} } {1 \over \ep + |0 - ( \e _{- i}^{(p)} )^{\top } \bbe |} \Big\} \prod_{i \in \Kc \setminus \{ i(1), \dots , i(p) \} } {1 \over \ep + | a_i - {\x _i}^{\top } \bbe |} \le M_5 {1 \over (1 + \| \bbe \| )^{| \Kc | - p + 1}} \non 
\end{align}
for all $\bbe \in \mathbb{R} ^p$ for some $M_5 > 0$. 
Note also that for some $M_6 > 0$, we have for all $\bbe \in \mathbb{R} ^p$ that $\| \bbe \| \ge \om / M_6$ whenever 
$1 - \prod_{i \in \Lc } 1(| y_i - {\x _i}^{\top } \bbe | \ge | y_i | / 2) = 1$. 
Then, since 
\begin{align}
f( y_i ) &\ge {1 \over M_7} {1 \over \om ( \log \om )^{1 + \ga }} \non 
\end{align}
for sufficiently large $\om > 0$ for all $i \in \Lc $ for some $M_7 > 0$ by part (iv) of Lemma \ref{lem:f}, 
\begin{align}
&h_{\i } ( \bbe , \si ; \om ) \non \\
&\le 1( \bbe \in A_{\i } ( \om )) \Big\{ 1 - \prod_{i \in \Lc } 1(| y_i - {\x _i}^{\top } \bbe | \ge | y_i | / 2) \Big\} \pi ( \si ) \Big\{ \prod_{j = 1}^{p} g_j ( \bbe , \si ; \om ) \Big\} M_8 {\om ^{| \Lc |} ( \log \om )^{| \Lc | (1 + \ga )} \over \om ^{| \Kc | - p + 1}} \non \\
&\le \pi ( \si ) \Big\{ \prod_{j = 1}^{p} g_j ( \bbe , \si ; \om ) \Big\} M_8 {\om ^{| \Lc |} ( \log \om )^{| \Lc | (1 + \ga )} \over \om ^{| \Kc | - p + 1}} \text{.} \non 
\end{align}
for sufficiently large $\om > 0$ and all $( \bbe , \si ) \in \mathbb{R} ^p \times (0, \infty )$ for some $M_8 > 0$. 
Thus, 
\begin{align}
\int_{\mathbb{R} ^p \times (0, \infty )} h_{\i } ( \bbe , \si ; \om ) d( \bbe , \si ) &\le M_8 \int_{\mathbb{R} ^p \times (0, \infty )} \pi ( \si ) \Big\{ \prod_{j = 1}^{p} g_j ( \bbe , \si ; \om ) \Big\} d( \bbe , \si ) {\om ^{| \Lc |} ( \log \om )^{| \Lc | (1 + \ga )} \over \om ^{| \Kc | - p + 1}} \non \\
&= M_9 {\om ^{| \Lc |} ( \log \om )^{| \Lc | (1 + \ga )} \over \om ^{| \Kc | - p + 1}} \to 0 \non 
\end{align}
as $\om \to \infty $ for some $M_9 > 0$. 
This completes the proof. 
\hfill$\Box$

\section{Appendix. 
Proof of Lemma \ref{lem:gagnon}}
Here, we prove Lemma \ref{lem:gagnon}. 

\bigskip

\noindent
{\bf Proof of Lemma \ref{lem:gagnon}.} \ \ Let $B_{i}^{(r)} ( \om ) = \{ \bbe \in \mathbb{R} ^p | | a_i + b_i \om - {\z _i}^{\top } \bbe | \le r \} $ for $\om \ge 0$ for $r > 0$ for $i = 1, \dots , m$. 
For $h \in \mathbb{N}$, let $I_h = \{ ( i_1 , \dots , i_h ) \in \{ 1, \dots , m \} ^h | \text{$i_j \neq i_{j'}$ for all $j, j' = 1, \dots , h$ with $j \neq j'$} \} $. 
For $( i_1 , \dots , i_p ) \in I_p$, let $f_{( i_1 , \dots , i_p )} \colon \mathbb{R} ^p \to \mathbb{R} ^p$ be the function defined by $f_{( i_1 , \dots , i_p )} ( \bbe ) = ( \z _{i_1} , \dots , \z _{i_p} )^{\top } \bbe $, $\bbe \in \mathbb{R} ^p$. 
Then for any $( i_1 , \dots , i_p ) \in I_p$, $f_{( i_1 , \dots , i_p )}$ is invertible by assumption (i) and, by linearity, 
\begin{align}
\bigcap_{j = 1}^{p} B_{i_j}^{(r)} ( \om ) &= {f_{( i_1 , \dots , i_p )}}^{- 1} (( a_{i_1} + b_{i_1} \om , \dots , a_{i_p} + b_{i_p} \om )^{\top } + [- r, r]^p ) \non \\
&= {f_{( i_1 , \dots , i_p )}}^{- 1} (( a_{i_1} + b_{i_1} \om , \dots , a_{i_p} + b_{i_p} \om )^{\top } ) +  {f_{( i_1 , \dots , i_p )}}^{- 1} ([- r, r]^p ) \label{lgagnonp1} 
\end{align}
for all $\om \ge 0$ for all $r > 0$. 

For all $( i_1 , \dots , i_p , i_{p + 1} ) \in I_{p + 1}$ and all $\bbe \in \mathbb{R} ^p$, we have by assumption (ii) that 
\begin{align}
&\begin{pmatrix} {\z _{i_1}}^{\top } \\ \vdots \\ {\z _{i_{p + 1}}}^{\top } \end{pmatrix} \bbe \neq \begin{pmatrix} a_{i_1} \\ \vdots \\ a_{i_{p + 1}} \end{pmatrix} \text{,} \non 
\end{align}
which implies that $| a_{i_j} - {\z _{i_j}}^{\top } \bbe | > \de $ for some $1 \le j \le p + 1$ for some $\de > 0$ and hence that 
\begin{align}
\bbe \in \bigcup_{k = 1}^{\infty } \bigcup_{j = 1}^{p + 1} \{ B_{i_j}^{(1 / k)} (0) \} ^c \text{.} \non 
\end{align}
Therefore, 
\begin{align}
\mathbb{R} ^p \subset \bigcap_{( i_1 , \dots , i_p , i_{p + 1} ) \in I_{p + 1}} \bigcup_{k = 1}^{\infty } \bigcup_{j = 1}^{p + 1} \{ B_{i_j}^{(1 / k)} (0) \} ^c \non 
\end{align}
and, since $I_{p + 1}$ is finite and since $B_{i}^{(r)} (0) \subset B_{i}^{( r' )} (0) $ for all $r' > r > 0$ for all $i = 1, \dots , m$, we have 
\begin{align}
\mathbb{R} ^p \subset \bigcup_{k = 1}^{\infty } \bigcap_{( i_1 , \dots , i_p , i_{p + 1} ) \in I_{p + 1}} \bigcup_{j = 1}^{p + 1} \{ B_{i_j}^{(1 / k)} (0) \} ^c \text{.} \label{lgagnonp2} %
\end{align}
Meanwhile, for all $( i_1 , \dots , i_p ) \in I_p$ and all $k \ge 1$, 
\begin{align}
\bigcap_{j = 1}^{p} B_{i_j}^{(1 / k)} (0) &= {f_{( i_1 , \dots , i_p )}}^{- 1} (( a_{i_1} , \dots , a_{i_p} )^{\top } + [- 1 / k, 1 / k]^p ) \non \\
&\subset {f_{( i_1 , \dots , i_p )}}^{- 1} (( a_{i_1} , \dots , a_{i_p} )^{\top } + [- 1, 1]^p ) \non 
\end{align}
by (\ref{lgagnonp1}) and the right-hand side is a compact set since $f_{( i_1 , \dots , i_p )}$ is a homeomorphism. 
Therefore, 
\begin{align}
\bigcup_{( i_1 , \dots , i_p ) \in I_p} \bigcup_{k = 1}^{\infty } \bigcap_{j = 1}^{p} B_{i_j}^{(1 / k)} (0) \subset B \text{,} \label{lgagnonp3} %
\end{align}
where 
\begin{align}
B = \bigcup_{( i_1 , \dots , i_p ) \in I_p} {f_{( i_1 , \dots , i_p )}}^{- 1} (( a_{i_1} , \dots , a_{i_p} )^{\top } + [- 1, 1]^p ) \text{.} \non 
\end{align}
Note that $B$ is compact since $I_p$ is finite. 
Then 
\begin{align}
B \subset \bigcap_{( i_1 , \dots , i_p , i_{p + 1} ) \in I_{p + 1}} \bigcup_{j = 1}^{p + 1} \{ B_{i_j}^{(1 / k_0 )} (0) \} ^c \non %
\end{align}
for some $k_0 \ge 1$ by (\ref{lgagnonp2}) and by the monotonicity of the sequence of open sets 
\begin{align}
\Big( \bigcap_{( i_1 , \dots , i_p , i_{p + 1} ) \in I_{p + 1}} \bigcup_{j = 1}^{p + 1} \{ B_{i_j}^{(1 / k)} (0) \} ^c \Big) _{k = 1}^{\infty } \text{.} \non 
\end{align}
However, by (\ref{lgagnonp3}), 
\begin{align}
B^c &\subset \bigcap_{( i_1 , \dots , i_p ) \in I_p} \bigcap_{k = 1}^{\infty } \bigcup_{j = 1}^{p} \{ B_{i_j}^{(1 / k)} (0) \} ^c \non \\
&\subset \bigcap_{( i_1 , \dots , i_p ) \in I_p} \bigcup_{j = 1}^{p} \{ B_{i_j}^{(1 / k_0 )} (0) \} ^c \non \\
&\subset \bigcap_{( i_1 , \dots , i_p , i_{p + 1} ) \in I_{p + 1}} \bigcup_{j = 1}^{p + 1} \{ B_{i_j}^{(1 / k_0 )} (0) \} ^c \text{.} \non %
\end{align}
Thus, 
\begin{align}
\mathbb{R} ^p \subset \bigcap_{( i_1 , \dots , i_p , i_{p + 1} ) \in I_{p + 1}} \bigcup_{j = 1}^{p + 1} \{ B_{i_j}^{(1 / k_0 )} (0) \} ^c \non 
\end{align}
and 
\begin{align}
\bigcup_{( i_1 , \dots , i_p , i_{p + 1} ) \in I_{p + 1}} \bigcap_{j = 1}^{p + 1} B_{i_j}^{( \ep )} (0) \subset \emptyset \non 
\end{align}
for some $\ep > 0$. 

Let $K = \{ i \in \{ 1, \dots , m \} | b_i = 0 \} $ and $L = \{ 1, \dots , m \} \setminus K = \{ i \in \{ 1, \dots , m \} | b_i \neq 0 \} $. 
Fix $( i_1 , \dots , i_p , i_{p + 1} ) \in I_{p + 1}$ and note that 
\begin{align}
\bigcap_{j = 1}^{p + 1} B_{i_j}^{( \ep )} (0) = \emptyset \text{.} \non 
\end{align}
First, suppose that $\{ i_1 , \dots , i_p , i_{p + 1} \} \subset K$. 
Then clearly 
\begin{align}
\bigcap_{j = 1}^{p + 1} B_{i_j}^{( \ep )} ( \om ) = \bigcap_{j = 1}^{p + 1} B_{i_j}^{( \ep )} (0) = \emptyset \non 
\end{align}
for all $\om \ge 0$. 
Next, suppose that $\{ i_1 , \dots , i_p , i_{p + 1} \} \not\subset K$ and $\{ i_1 , \dots , i_p , i_{p + 1} \} \not\subset L$. 
Then there exist $1 \le j_1 , j_2 \le p + 1$ such that $i_{j_1} \in K$, $i_{j_2} \in L$, and $j_1 \neq j_2$. 
Let $( {i_1}' , \dots , {i_p}' ) \in I_p$ satisfy $\{ {i_1}' , \dots , {i_p}' \} = \{ i_1 , \dots , i_p , i_{p + 1} \} \setminus \{ i_{j_1} \} $ and let $\Zbt = ( \z _{{i_1}'} , \dots , \z _{{i_p}'} )^{\top }$, $\abt = ( a_{{i_1}'} , \dots , a_{{i_p}'} )^{\top }$, and $\bbt = ( b_{{i_1}'} , \dots , b_{{i_p}'} )^{\top }$. 
Then, by (\ref{lgagnonp1}), 
\begin{align}
\bigcap_{j = 1}^{p + 1} B_{i_j}^{( \ep )} ( \om ) &= B_{i_{j_1}}^{( \ep )} ( \om ) \cap \bigcap_{j = 1}^{p} B_{{i_j}'}^{( \ep )} ( \om ) \non \\
&= \{ \bbe \in \mathbb{R} ^p | | a_{i_{j_1}} + b_{i_{j_1}} \om - {\z _{i_{j_1}}}^{\top } \bbe | \le \ep \} \non \\
&\quad \cap \{ {f_{( {i_1}' , \dots , {i_p}' )}}^{- 1} (( a_{{i_1}'} + b_{{i_1}'} \om , \dots , a_{{i_p}'} + b_{{i_p}'} \om )^{\top } ) + {f_{( {i_1}' , \dots , {i_p}' )}}^{- 1} ([- \ep , \ep ]^p ) \} \non \\
&= \{ \bbe \in \mathbb{R} ^p | | a_{i_{j_1}} - {\z _{i_{j_1}}}^{\top } \bbe | \le \ep \} \cap \{ ( \Zbt ^{- 1} \abt + \Zbt ^{- 1} \bbt \om ) + {f_{( {i_1}' , \dots , {i_p}' )}}^{- 1} ([- \ep , \ep ]^p ) \} \non 
\end{align}
for all $\om \ge 0$. 
Note that $\Zbt $ is invertible by assumption (i) and that 
\begin{align}
\bm{0} ^{(p + 1)} \neq \begin{pmatrix} \Zbt & \bbt \\ {\z _{i_{j_1}}}^{\top } & b_{i_{j_1}} \end{pmatrix} \begin{pmatrix} \Zbt ^{- 1} \bbt \\ - 1 \end{pmatrix} = \begin{pmatrix} \Zbt & \bbt \\ {\z _{i_{j_1}}}^{\top } & 0 \end{pmatrix} \begin{pmatrix} \Zbt ^{- 1} \bbt \\ - 1 \end{pmatrix} = \begin{pmatrix} \bm{0} _p \\ {\z _{i_{j_1}}}^{\top } \Zbt ^{- 1} \bbt \end{pmatrix} \non 
\end{align}
by assumption (iii). 
Then ${\z _{i_{j_1}}}^{\top } \Zbt ^{- 1} \bbt \neq 0$ 
and it follows that for any $\om \ge 0$, we have that for all $\bbe \in \bigcap_{j = 1}^{p + 1} B_{i_j}^{( \ep )} ( \om )$, there exists $\v \in {f_{( {i_1}' , \dots , {i_p}' )}}^{- 1} ([- \ep , \ep ]^p )$ such that 
\begin{align}
\ep &\ge | a_{i_{j_1}} - {\z _{i_{j_1}}}^{\top } \bbe | = | a_{i_{j_1}} - {\z _{i_{j_1}}}^{\top } ( \Zbt ^{- 1} \abt + \Zbt ^{- 1} \bbt \om + \v )| \non \\
&\ge | {\z _{i_{j_1}}}^{\top } \Zbt ^{- 1} \bbt | \om -  | a_{i_{j_1}} | - | {\z _{i_{j_1}}}^{\top } \Zbt ^{- 1} \abt | - \| {\z _{i_{j_1}}} \| \| \v \| \text{.} \non 
\end{align}
Thus, there exists $M > 0$ such that $\bigcap_{j = 1}^{p + 1} B_{i_j}^{( \ep )} ( \om ) = \emptyset $ for all $\om \ge M$. 
Finally, suppose that $\{ i_1 , \dots , i_p , i_{p + 1} \} \not\subset K$ and $\{ i_1 , \dots , i_p , i_{p + 1} \} \subset L$. 
Let $( {i_1}' , \dots , {i_p}' ) = ( i_1 , \dots , i_p ) \in I_p$ and let $\Zbt = ( \z _{{i_1}'} , \dots , \z _{{i_p}'} )^{\top }$, $\abt = ( a_{{i_1}'} , \dots , a_{{i_p}'} )^{\top }$, and $\bbt = ( b_{{i_1}'} , \dots , b_{{i_p}'} )^{\top }$. 
Similarly, let $( {i_1}'' , \dots , {i_p}'' ) = ( i_1 , \dots , i_{p - 1} , i_{p + 1} ) \in I_p$ and let $\widetilde{\Zbt } = ( \z _{{i_1}''} , \dots , \z _{{i_p}''} )^{\top }$, $\tilde{\abt } = ( a_{{i_1}''} , \dots , a_{{i_p}''} )^{\top }$, and $\tilde{\bbt } = ( b_{{i_1}''} , \dots , b_{{i_p}''} )^{\top }$. 
Then, by (\ref{lgagnonp1}), 
\begin{align}
\bigcap_{j = 1}^{p + 1} B_{i_j}^{( \ep )} ( \om ) &= \{ ( \Zbt ^{- 1} \abt + \Zbt ^{- 1} \bbt \om ) + {f_{( {i_1}' , \dots , {i_p}' )}}^{- 1} ([- \ep , \ep ]^p ) \} \non \\
&\quad \cap \{ ( \widetilde{\Zbt } ^{- 1} \tilde{\abt } + \widetilde{\Zbt } ^{- 1} \tilde{\bbt } \om ) + {f_{( {i_1}'' , \dots , {i_p}'' )}}^{- 1} ([- \ep , \ep ]^p ) \} \non 
\end{align}
for all $\om \ge 0$. 
Note that $\Zbt $ and $\widetilde{\Zbt }$ are invertible by assumption (i) and that 
\begin{align}
&\begin{pmatrix} {\z _{{i_p}''}}^{\top } & b_{{i_p}''} \end{pmatrix} \neq {\z _{{i_p}''}}^{\top } \Zbt ^{- 1} \begin{pmatrix} \Zbt & \bbt \end{pmatrix} \non 
\end{align}
by assumption (iii). 
Then 
\begin{align}
\begin{pmatrix} \widetilde{\Zbt } & \tilde{\bbt } \end{pmatrix} \neq \widetilde{\Xbt } \Zbt ^{- 1} \begin{pmatrix} \Zbt & \bbt \end{pmatrix} \quad \text{and} \quad \widetilde{\Zbt } ^{- 1} \tilde{\bbt } \neq \Zbt ^{- 1} \bbt \text{.} \non 
\end{align}
Therefore, for any $\om \ge 0$, we have that for all $\bbe \in \bigcap_{j = 1}^{p + 1} B_{i_j}^{( \ep )} ( \om )$, there exist $\vbt \in {f_{( {i_1}' , \dots , {i_p}' )}}^{- 1} ([- \ep , \ep ]^p )$ and $\tilde{\vbt } \in {f_{( {i_1}'' , \dots , {i_p}'' )}}^{- 1} ([- \ep , \ep ]^p )$ such that $\Zbt ^{- 1} \abt + \Zbt ^{- 1} \bbt \om + \vbt = \bbe = \widetilde{\Zbt } ^{- 1} \tilde{\abt } + \widetilde{\Zbt } ^{- 1} \tilde{\bbt } \om + \tilde{\vbt }$, which implies that 
$( \Zbt ^{- 1} \bbt - \widetilde{\Zbt } ^{- 1} \tilde{\bbt } ) \om = \widetilde{\Zbt } ^{- 1} \tilde{\abt } - \Zbt ^{- 1} \abt + \tilde{\vbt } - \vbt$ and hence that $\om \le ( \| \widetilde{\Zbt } ^{- 1} \tilde{\abt } - \Zbt ^{- 1} \abt \| + \| \tilde{\vbt } \| + \| \vbt \| ) / \| \Zbt ^{- 1} \bbt - \widetilde{\Zbt } ^{- 1} \tilde{\bbt } \| $. 
Thus, there exists $M > 0$ such that $\bigcap_{j = 1}^{p + 1} B_{i_j}^{( \ep )} ( \om ) = \emptyset $ for all $\om \ge M$. 

Since $I_{p + 1}$ is finite, we conclude that there exists $M > 0$ such that for all $\om \ge M$, 
\begin{align}
\bigcup_{( i_1 , \dots , i_p , i_{p + 1} ) \in I_{p + 1}} \bigcap_{j = 1}^{p + 1} B_{i_j}^{( \ep )} ( \om ) = \emptyset \text{.} \non 
\end{align}
Hence, for all $\om \ge M$, 
\begin{align}
\mathbb{R} ^p &= \bigcap_{( i_1 , \dots , i_p , i_{p + 1} ) \in I_{p + 1}} \bigcup_{j = 1}^{p + 1} \{ B_{i_j}^{( \ep )} ( \om ) \} ^c \non \\
&\subset \bigcup_{i_1 = 1}^{m} \dots \bigcup_{i_p = 1}^{m} \bigcap_{i \in \{ 1, \dots , m \} \setminus \{ i_1 , \dots , i_p \} } \{ \bbe \in \mathbb{R} ^p | | a_i + b_i \om - {\x _i}^{\top } \bbe | > \ep \} \text{.} \non 
\end{align}
This completes the proof. 
\hfill$\Box$

\section*{Acknowledgments}
Research of the author was supported in part by JSPS KAKENHI Grant Number JP22K20132, JP19K11852 from Japan Society for the Promotion of Science.

\newpage
\setcounter{page}{1}
\setcounter{equation}{0}
\renewcommand{\theequation}{S\arabic{equation}}
\setcounter{section}{0}
\renewcommand{\thesection}{S\arabic{section}}
\setcounter{table}{0}
\renewcommand{\thetable}{S\arabic{table}}
\setcounter{figure}{0}
\renewcommand{\thefigure}{S\arabic{figure}}

\begin{center}
{\LARGE\bf Supplementary Materials}
\end{center}

\bigskip

\section{Properties of the Error Density of Section \ref{sec:conditions} of the Main Text}
The following result is used in Section \ref{sec:proof} of the main text.

\begin{lem}
\label{lem:f} 
Let $\ga > 0$ and let 
\begin{align}
f(z) = {\ga / 2 \over 1 + |z|} {1 \over \{ 1 + \log (1 + |z|) \} ^{1 + \ga }} \non 
\end{align}
for $z \in \mathbb{R}$. 
\begin{itemize}
\item[{\rm{(i)}}]
For all $( \mu , \si ) \in \mathbb{R} \times (0, \infty )$, 
\begin{align}
f((y - \mu ) / \si ) / \si \sim f(y) \non 
\end{align}
as $y \to \pm \infty $. 
\item[{\rm{(ii)}}]
Let $( \mu , \si ) \in \mathbb{R} \times (0, \infty )$ and $y \in \mathbb{R}$. 
Then, if $|y - \mu | \ge |y| / 2$ and if $|y| \ge 1$, we have 
\begin{align}
{f((y - \mu ) / \si ) / \si \over f(y)} &\le 4 (1 + \log 3)^{1 + \ga } \{ 1 + \log (1 + \si ) \} ^{1 + \ga } \text{.} \non 
\end{align}
\item[{\rm{(iii)}}]
For all $( \mu , \si ) \in \mathbb{R} \times (0, \infty )$ and all $y \in \mathbb{R}$, 
\begin{align}
f((y - \mu ) / \si ) / \si \le {\ga \over 2} {1 \over \si } {1 \over 1 + |y - \mu | / \si } \text{.} \non 
\end{align}
\item[{\rm{(iv)}}]
For all $y \in \mathbb{R}$ with $|y| \ge 2 e$, 
\begin{align}
f(y) %
&\ge {\ga \over 2^{3 + \ga }} {1 \over |y| ( \log |y|)^{1 + \ga }} \text{.} \non 
\end{align}
\end{itemize}
\end{lem}

\noindent
{\bf Proof%
.} \ \ Part (i) follows from Lemma S1 of Hamura et al. (2022). 
Part (ii) follows since $\{ 1 + \log (1 + y_1 y_2 ) \} / \{ 1 + \log (1 + y_1 ) \} \le 1 + \log (1 + y_2 )$ for all $y_1 , y_2 > 0$; this inequality is proved in the proof of Lemma S6 of Hamura et al. (2022). 
Parts (iii) and (iv) are trivial. 
\hfill$\Box$

\end{document}